*Original Article*

# Reconsideration of Tangle and Ultrafilter using Separation and Partition

Takaaki Fujita

*Graduate school of Science and Technology, Gunma University, 1-5-1 Tenjin-cho Kiryu Gunma, Japan.*

*Corresponding Author : t171d603@gunma-u.ac.jp*

***Abstract -*** *Tangle is a concept in graph theory that has a dual relationship with branch-width which is well-known graph width parameter. Ultrafilter, a fundamental notion in mathematics, is similarly known to have a dual relationship with branch-width when extended to a connectivity system (X, f). We will reconsider these concepts using separation and partition.*

## 1. Introduction

Graph theory, a fundamental branch of mathematics, is focused on analyzing networks composed of nodes and edges, exploring their paths, structures, and properties. Within this field, a critical metric known as the "graph width parameter" measures the width of a graph. This parameter represents the maximum width across all cuts or layers in a hierarchical decomposition of the graph, playing a pivotal role in assessing the complexity and structural properties of graphs. The significance of graph width parameters extends across various disciplines, such as computer science, network theory, artificial intelligence, matroid theory, and discrete mathematics. These parameters, including branch-width, tree-width, carving-width, path-distance-width, linear-width, and path-width, are extensively documented in graph theory and combinatorial literature. The vast body of research surrounding these metrics underscores their importance and broad applications, as indicated by numerous studies [4, 5, 6, 10, 11, 12, 13, 14, 15, 18, 19, 20, 21, 24, 28, 30, 34, 54, 55]. In practical applications, particularly in engineering and computer science, graph width parameters such as tree-width and branch-width are crucial. The concept of tree-width, for example, is instrumental in simplifying complex network analysis by enabling the breakdown of intricate graphs into more manageable, tree-like structures. This decomposition significantly benefits algorithms that would otherwise be inefficient on general graphs. Algorithms that operate within polynomial time on graphs with bounded tree-width, for instance, can efficiently solve NP-hard problems such as the Hamiltonian cycle, maximum independent set, and graph coloring. Additionally, in the study of graph width parameters, concepts like duality theorems and equivalence relations are frequently explored. A duality theorem posits that the existence (or non-existence) of one entity implies the non-existence (or existence) of its counterpart. Similarly, equivalence relations determine that the existence (or non-existence) of one entity ensures the existence (or non-existence) of another specific entity. These theoretical frameworks are crucial for demonstrating the lower bounds of width parameters in minimization problems, providing valuable insights into the structural properties and optimization potential of graphs. Tangle, a concept first introduced in reference [1], plays an essential role in determining whether a branch-width is at most k (k is natural number), where $k+1$ represents the order of the tangle (see also references [2, 3, 33]). Tangle has a crucial role in graph theory and has been extensively studied by many mathematicians, making its research inherently important (ex. [4, 7, 25, 26]). In this concise paper, a pair $(X, f)$ of a finite set (an underlying set) $X$ and a symmetric submodular function $f$ is called a connectivity system. The ultrafilter, a well-known mathematical notion, is also known to have a dual relationship with branch-width when extended over a connectivity system $(X, f)$.

Tangle can be more comprehensibly and intuitively interpreted using set theory concepts such as separation and partition. In literature [16, 17, 27, 28], the authors investigate the profound relationship between width parameters and Tangle (Profile) using the concept of an abstract separation system. By using separation and partition, the states on both the left and right sides become intuitively easier to understand, not just on one side. In this concise paper, we will reconsider the relationship between Tangle and Ultrafilter using partition and separation. Although it has been demonstrated that Tangle and Ultrafilter have a cryptomorphism on a connectivity system $(X, f)$, we will reconsider this relationship using partition and separation. Also we discussed about weak ultrafilter and profile using partition and separation on a connectivity system $(X, f)$. At first glance, the idea that seemingly unrelated concepts can become relevant when considered within the context of a connectivity system is both profound and surprising. By reexamining the dual concepts of graph width parameters, this research aims to contribute to the study of graph width parameters and graph algorithms.

## 2. Preliminaries

In this section, we compile the necessary definitions for this paper. In this paper, we use expressions like $A \subseteq X$ to indicate that A is a subset of $X$, $A \cup B$ to represent the union of two subsets $A$ and $B$, both of which are subsets of $X$, or $A = \emptyset$ to signify an empty set. Specifically, $A \cap B$ denotes the intersection of subsets $A$ and $B$. A similar logic applies to $A \setminus B$.

### *2.1. Symmetric Submodular Function*

Symmetric submodular functions constitute an essential family of submodular functions, encompassing a range of fascinating cases, including hypergraphs and graph width parameters. Numerous studies have also been conducted on symmetric submodularity (e.g., [22, 23]). The definition of a symmetric submodular function is presented below.

**Definition 1:** Let $X$ be a finite set. A function $f: X \rightarrow \mathbb{N}$ is called symmetric submodular if it satisfies the following conditions: · $\forall A \subseteq X, f(A) = f(X \setminus A)$ (symmetry), · $\forall A, B \subseteq X, f(A) + f(B) \geq f(A \cap B) + f(A \cup B)$ (submodularity).

In this short paper, a pair $(X, f)$ of a finite set (an underlying set) $X$ and a symmetric submodular function $f$ is called a connectivity system. In this paper, we use the notation $f$ for a symmetric submodular function, a finite set (referred to as the underlying set) $X$, and a natural number $k$. A set $A$ is $k$-efficient if $f(A) \leq k$. A symmetric submodular function possesses the following principalties. This lemma will be utilized in the proofs of lemmas and theorems presented in this paper.

**Lemma 2 [10]:** A symmetric submodular function f satisfies:
1. $\forall A \subseteq X, f(A) \geq f(\emptyset) = f(X)$,
2. $\forall A, B \subseteq X, f(A) + f(B) \geq f(A \setminus B) + f(B \setminus A)$.

### *2.2. Partitions and separations of sets*

In terms of graph width parameters, concepts such as partitions and separation are sometimes used to discuss these parameters in a more understandable and abstract way. This paper also considers graph width parameters using the concepts of partitions and separation. A partition refers to the process of dividing the elements of a set into non-empty, distinct subsets such that each element belongs to one and only one subset. Within this treatise, a partition $(A, B)$ in a pair $(X, f)$ is designated as a separation in a pair $(X, f)$ (cf. [2]), with $A$ and $B$ referred to as sides of said separation. An order of a separation $(A, B)$ in a pair $(X, f)$ is delineated as $f(A)$, wherein $f$ represents a symmetric submodular function. Every separation $(A, B)$ comprising subsets of a finite set $X$ adheres to $A \cup B = X$. We establish a partial order denoted by $\leq$ on the set of separations of an underlying set $X$ as follows:

**Notation 3 (cf. [1]):** Given two separations $(A, B)$ and $(C, D)$ of a finite set $X$, we permit $(A, B) \leq (C, D)$ provided that $A \subseteq C$ and $B \supseteq D$. We inscribe $(A, B) < (C, D)$ exclusively when $(A, B) \leq (C, D)$ and $(A, B) \neq (C, D)$.

In accordance with this characterization, we deduce that $(A, B) \leq (C, D)$ if and only if $(B, A) \geq (D, C)$. The above definition of partial order employs the approach from reference [1]; however, it differs in that the separation in reference [1] is not limited to partition. In this paper, we would not consider about an abstract separation system.

### *2.3. Tangle of separations*

The tangle serves as a dual concept to branch-decomposition, a well-known width parameter of a graph [33]. Tangles can be defined not only on graphs but also for symmetric submodular functions $f$ over a finite set $X$ mapping to non-negative integers [35]. There are several combinatorial definitions deeply related to tangles. The definition of a tangle of separations is provided below:

**Definition 4 [2]:** A tangle of separations in $(X,f)$ of order $k + 1$ is a family $T$ of separations of finite set $X$, each of order $\leq k$, such that:
(T1) For every separation $(A, B)$ of $X$ of order $\leq k$, one of $(A, B)$, $(B, A)$ is an element of $T$;
(T2) For each element $e \in X$ such that $f(\{e\}) \leq k$ , $(e, E \setminus \{e\}) \in T$;
(T3) If $(A_1, B_1), (A_2, B_2), (A_3, B_3) \in$ T then $A_1 \cup A_2 \cup A_3 \neq X$.

Related to the above concept, Linear Tangle is known as a concept with a dual relationship to Linear Width [3, 20]. A Linear Tangle of separation can be obtained by transforming the axiom (T3) of the Tangle definition mentioned above to the following axiom (LT3). This axiom (LT3) is called single-element-extension axiom (also lifting axiom) in this paper (cf. [31, 32, 40]).
(LT3) $(A_1, B_1), (A_2, B_2) \in T, e \in X, f(\{e\}) \leq k \Rightarrow A_1 \cup A_2 \cup \{e\} \neq X.$

In referenc3e [2], where the Tangle Bases is defined, a Linear Tangle Bases can be established by modifying the axiom (T3) to axiom (LT3), thereby limiting it to singletons.

A tangle and a linear tangle have the following principalties.

**Lemma 5:** A tangle T of separations in (X,f) of order k + 1 satisfies the following axiom (T4).
(T4) $(\emptyset, X) \in T$
**Proof of lemma 5:** We prove this by contradiction. Assume that $(\emptyset, X) \notin T$. Then by Lemma 2 and the axiom (T1), we have $(X, \emptyset) \in T$, which contradicts axiom (T3). Hence, axiom (T4) must hold. This proof is completed.

### *2.4. Ultrafilter of separations*

An ultrafilter on a set is a maximal filter in a partially ordered set, typically used in mathematical contexts like topology and logic. It selects precisely one subset from each partition of a set, thus serving as a decisive tool for examining properties of sets and functions. We will define the Ultrafilter of separations on $(X, f)$ using separations. This new definition maintains a dual relationship with branch-width and holds an equivalent relationship with Tangle (see section 4). The definition of a separation ultrafilter, interpreted through separation, is provided below by adding the condition of a submodular function to ultrafilters on a Boolean algebra.

**Definition 6:** Let $X$ represent a finite set and $f$ denote a symmetric submodular function delineated over $X$, mapping to the non-negative integers. An Ultrafilter of separations in $(X,f)$ of order $k + 1$ is a family $F$ of separations of finite set $X$, each of order $\leq k$, such that:
(F1) For every separation $(A, B)$ of $X$ of order $\leq k$ , one of $(A, B), (B, A)$ is an element of $F$,
(F2) $(\emptyset, X) \notin F$,
(F3) For each element $e \in X$ such that $f(\{e\}) \leq k$ , $(e, E \setminus \{e\}) \notin F$;
(F4) If $(A_1, B_1) \in$ F, $(A_1, B_1) \leqq (A_2, B_2)$ , order of $(A_2, B_2) \leq k$ then $(A_2, B_2) \in$ F,
(F5) If $(A_1, B_1) \in F, (A_2, B_2) \in F$, order of $(A_1 \cap A_2, B_1 \cup B_2) \leq k$ then $(A_1 \cap A_2, B_1 \cup B_2) \in F,$

By modifying axiom (F5) in the above Ultrafilter definition to the following (SF5), we can define a Single UltraFilter of separations, which has a dual relationship with Linear-width and holds an equivalent relationship with Linear Tangle (see section 4). This axiom is called single-element-deletion axiom in this paper. And this concept is deep relation to single ideal (cf. [3,21,38]).
(SF5) If $(A_1, B_1) \in F, e \in X, f(\{e\}) \leq k,$ order of $(A_1 \cap (X \setminus \{e\}), B_1 \cup \{e\}) \leq k$ *then* $(A_1 \cap (X \setminus \{e\}), B_1 \cup \{e\}) \in F,$

Ultrafilters have a wide range of engineering applications, which has attracted the attention of numerous researchers (ex. [46,47,48]). Therefore, it is believed that studying ultrafilters holds significant value. Note that Ultrafilter on $(X, f)$ is co-maximal ideal on $(X, f)$.

## 3. Consideration of Bramble and Weak Ultrafilter of separations

This section discuss about Bramble and Weak Ultrafilter of separations. One of well-known concept in graph width parameters is called Bramble, which has a close relationship with Tangle. Intuitively speaking, fot every separation $(A, B)$ of $X$ of order $\leq k$, then the tangle captures the "small" part of separation $(A, X \setminus A)$ while the bramble encompasses the "large (big)" one. Due to the relationship between Tangle and Bramble, it can be said that Bramble is a dual concept to the width parameter such as branch-width, tree-width and so on. Numerous studies have been conducted on this particular Bramble (ex.[41,42,43,44,45]). An Ultrafilter of separations satisfies certain conditions.

**Theorem 7:** Let $X$ represent a finite set and $f$ denote a symmetric submodular function. *If a family F* is an Ultrafilter of separations of order $k+1$ on $(X,f)$, We obtain following axiom.

(F6) If $(A_1, B_1) \in F$, $(A_2, B_2) \in F$, $(A_3, B_3) \in F$ then $A_1 \cap A_2 \cap A_3$ is non-empty.

**Proof of Theorem 7:** Suppose, for the sake of contradiction, that there exist $(A_1, B_1)$, $(A_2, B_2)$, and $(A_3, B_3) \in$ F with $A_1 \cap A_2 \cap A_3 = \emptyset$. Choose $A_1 \cap A_2 \cap A_3$ to be inclusion-wise maximal with $f(A) \leq k$. We now claim that $f(A_1 \cap A_3) \leq k$ *and* $f(A_1 \setminus A_3) \leq k$. First, we show that $f(A_1 \cap A_3) \leq k$. Indeed, let $(C, D) = ((A_1 \setminus A_2) \cup A_3, X \setminus ((A_1 \setminus A_2) \cup A_3))$. Suppose $(A_3, B_3) \leq (C, D)$, meaning $(A_1 \setminus A_2, X \setminus (A_1 \setminus A_2)) \leq (C, D)$. Since $A_1 \cap A_2 \cap A_3 = \emptyset$, we obtain $A_1 \setminus A_2 = A_1 \cap A_3$. So we obtain $f(A_1 \cap A_3) \leq k$ and the claim holds. Suppose $(A_3, B_3) \not\leq (C, D)$. If $f(C) \leq k$, then $(C, D) \in F$. However, this is impossible since $A_1$, $A_2$, and $A_3$ are maximal with this principalty. Thus, we consider $f(C) > k$. And we have $f(C) + f((A_1 \setminus A_2) \cap A_3) \leq f(A_1 \setminus A_2) + f(A_3) \leq 2k$. Therefore, we obtain $f((A_1 \setminus A_2) \cap A_3) \leq$ k. Finally, since $A_1 \cap A_2 \cap A_3 = \emptyset$, *we obtain* $A_1 \cap A_3 = (A_1 \setminus A_2) \cap A_3$. So we obtain $f(A_1 \cap A_3) \leq k$ and the claim holds. Using the same calculations as previously discussed, it can be assumed that $f(A_1 \setminus A_3) \leq k$. However, the separation(partition) (X\A, $A_1 \cap A_3$, $A_1 \setminus A_3$) is not possible, as each of these three sets is disjoint from $A_1$, $A_2$, and $A_3$, which are all members of $A_2$. Therefore, we have demonstrated that if $(A_1, B_1) \in F$, $(A_2, B_2) \in F$, $(A_3, B_3) \in F$, then $A_1 \cap A_2 \cap A_3$ is non-empty. The proof is now complete. By modifying axiom (F5) in the Ultrafilter definition above to the following axiom (WF5), we can define a Weak Ultrafilter of separations. The Weak Ultrafilter is a concept used in the world of logic [36, 37, 39]. It is worth noting that, in general, the complement set of a Weak filter is known to be a Weak ideal, and this can be similarly defined on $(X, f)$ as well. (WF5) If $(A_1, B_1) \in F$, $(A_2, B_2) \in F$, order of $(A_1 \cap A_2, B_1 \cup B_2) \leq k$, then $A_1 \cap A_2 \neq \emptyset$. The following Theorem 8 clearly holds true.

**Theorem 8:** Let $X$ represent a finite set and $f$ denote a symmetric submodular function. If W is an Ultrafilter of separations of order $k+1$, *then W is* a Weak Ultrafilter of separations of order $k+1$.

**Proof of Theorem 8:** Since $W$ is an Ultrafilter of separations of order $k+1$, it already satisfies axiom (F1) to (F4). We only need to show that it satisfies axiom (WF5). Consider any two separations $(A_1, B_1)$ and $(A_2, B_2)$ in $W$. Now, if the order of $(A_1 \cap A_2, B_1 \cup B_2) \leq$ k, then by axiom (F2) and (F5), $A_1 \cap A_2$ must be non-empty. This means that W satisfies axiom (WF5). The proof is now complete.

Due to the principalties mentioned above, the author believes that there is a discomforting relationship between Weak Ultrafilters and graph width parameters such as Tangle, Bramble, and Branch Decomposition. However, a proof has not yet been established. The following is assumed to hold for linear tangles / single ultrafilter as well, naturally.

**Open Problem 9:** Let $X$ represent a finite set and $f$ denote a symmetric submodular function. *If a family W* is an Weak Ultrafilter of separations of order $k+1$ on $(X,f)$, We obtain following axiom.

(F6) If $(A_1, B_1) \in W$, $(A_2, B_2) \in W$, $(A_3, B_3) \in W$ then $A_1 \cap A_2 \cap A_3$ is non-empty.

**Open Problem 10:** Let $X$ represent a finite set and $f$ denote a symmetric submodular function. $T$ is a Tangle of separations of order $k+1$ on $(X,f)$ iff $W = \{(A,B) \mid (B,A) \in T\}$ is an Weak Ultrafilter of separations of order $k+1$ on $(X,f)$. Ultimately, our future goal is to demonstrate that a weak ultrafilter exists if and only if a non-principal bramble exists.

## 4. Cryptomorphism between Tangle of separations and Ultrafilter of separations

In this section, we demonstrate the cryptomorphism between Tangle of separations and Ultrafilter of separations. The main result of this paper is presented below. At first glance, the idea that seemingly unrelated concepts can become relevant when considered within the context of a connectivity system is both profound and surprising.

**Theorem 11:** Let $X$ represent a finite set and $f$ denote a symmetric submodular function. $T$ is a Tangle of separations of order $k+1$ on $(X,f)$ iff $F = \{(A,B) \mid (B,A) \in T\}$ is an Ultrafilter of separations of order $k+1$ on $(X,f)$.

**Proof of Theorem 11:** To prove this theorem, we need to show that if $T$ is a Tangle of separations of order $k+1$, then $F$ is an Ultrafilter of separations of order $k+1$, and vice versa.

($\Rightarrow$) Suppose $T$ is a Tangle of separations of order $k+1$ on $(X, f)$. By Definition 3, we know that $T$ satisfies axioms (T1), (T2), and (T3). We need to show that $F$ also satisfies the axioms (F1) - (F5) for Ultrafilters from Definition 4. We show axiom (F1) for Ultrafilters. From axiom (T1), we know that for every separation $(A, B)$ of $X$ of order $\leq k$, one of $(A, B)$, $(B, A)$ is an element of $T$. Since $F = \{(A,B) \mid (B,A) \in T\}$, axiom (F1) holds for $F$. We show axiom (F2) for Ultrafilters. The empty set $\emptyset$ cannot be an element of $F$ because the separations in $T$ have empty set $\emptyset$ by Lemma 4.

Thus, the axiom (F2) holds for $F$. We show axiom (F3) for Ultrafilters. Since axiom (T2) states that for each element $e \in X$ such that $f(\{e\}) \leq k$, $(e, X \setminus \{e\}) \in T$, it follows that $(X \setminus \{e\}, e) \notin F$. Thus, axiom (F3) holds for F. We now prove the Ultrafilter axiom (F4). Suppose $(A_1, B_1) \in F$ and $(A_1, B_1) \leq (A_2, B_2)$ with the order of $(A_2, B_2) \leq k$ and $(A_2, B_2) \notin F$. By assumption, $(B_1, A_1) \in T$ and $(A_2, B_2) \in T$. However, $B_1 \cup A_2 = X$, which contradicts the Tangle of separations axiom (T3). Therefore, the axiom (F4) is satisfied. We show the axiom (F5) for Ultrafilters. Suppose $(A_1, B_1) \in F$ and $(A_2, B_2) \in F$ with the order of $(A_1 \cap A_2, B_1 \cup B_2) \leq k$ and $(A_1 \cap A_2, B_1 \cup B_2) \notin F$. Then, $(B_1, A_1) \in T$ and $(B_2, A_2) \in T$. Moreover, by assumption, $(A_1 \cap A_2, B_1 \cup B_2) \in T$. By axiom (T3), $B_1 \cup B_2 \cup (A_1 \cap A_2) \neq X$. However, $B_1 \cup B_2 \cup (A_1 \cap A_2) = X$, which is a contradiction. Thus, $(A_1 \cap A_2, B_1 \cup B_2) \in F$, and axiom (F5) holds for $F$. Since F satisfies axioms (F1) - (F5), $F$ is an Ultrafilter of separations of order $k+1$ on $(X, f)$. (⇐) Conversely, suppose $F$ is an Ultrafilter of separations of order $k+1$ on $(X, f)$. Define $T = \{(A, B) \mid (B, A) \in F\}$. We need to show that $T$ satisfies the axioms (T1), (T2), and (T3) for Tangles from Definition 3. To show that $T$ satisfies axiom (T1), consider any separation $(A, B)$ of $X$ of order $\leq k$. By axiom (F1), we know that either $(A, B) \in F$ or $(B, A) \in F$. If $(A, B) \in F$, then $(B, A) \in T$, and if $(B, A) \in F$, then $(A, B) \in T$. Thus, $T$ satisfies axiom (T1). To show that T satisfies axiom (T2), consider any element $e \in X$ such that $f(\{e\}) \leq k$. By axiom (F3), we know that $(e, X \setminus \{e\}) \notin F$. Since either $(e, X \setminus \{e\}) \in F$ or $(X \setminus \{e\}, e) \in F$ by axiom (F1), it must be that $(X \setminus \{e\}, e) \in F$. Then, $(e, X \setminus \{e\}) \in T$. Thus, $T$ satisfies axiom (T2). To demonstrate that $T$ satisfies (T3), assume for the sake of contradiction that there exist separations $(A_1, B_1)$, $(A_2, B_2)$, and $(A_3, B_3)$ in T such that $A_1 \cup A_2 \cup A_3 = X$. This implies that $(B_1, A_1)$, $(B_2, A_2)$, and $(B_3, A_3)$ are in $F$. By employing the principalties of symmetric submodularity, the order of $(B_1 \cap B_2, A_1 \cup A_2)$ is equal to the order of $(A_3, B_3)$. Consequently, using the axiom (F5) of Ultrafilters, we obtain $(B_1 \cap B_2, A_1 \cup A_2) \in F$. Now, since $(B_1 \cap B_2, A_1 \cup A_2) \in F$ and $(B_3, A_3) \in F$, and the order of $((B_1 \cap B_2) \cap B_3, (B_1 \cap B_2, A_1 \cup A_2) = (\emptyset, X) \leq k$ by lemma 2, we can apply the axiom (F5) again to deduce that $((B_1 \cap B_2) \cap B_3, (A_1 \cup A_2) \cup A_3) \in F$. However, $(B_1 \cap B_2) \cap B_3 = \emptyset \in F$. This contradicts axiom (F2) for Ultrafilters, which states that the empty set $\emptyset$ cannot be an element of $F$. Thus, $T$ satisfies the axiom (T3). Since $T$ satisfies axioms (T1), (T2), and (T3), $T$ is a Tangle of separations of order $k+1$ on $(X, f)$. This proof is completed.From the aforementioned Theorem 11, Theorem 12 also holds true.

**Theorem 12:** Let $X$ represent a finite set and $f$ denote a symmetric submodular function delineated over $X$. $T$ is a Linear Tangle of separations of order $k+1$ on $(X,f)$ iff $F = \{(A,B) \mid (B,A) \in T\}$ is an Single Ultrafilter of separations of order $k+1$ on $(X,f)$.

## 5. Cryptomorphism between Profile of separations and Ultrafilter of separations

In this section, we also consider the extension of the concept of a Profile, initially defined for graphs in reference [49], to the connectivity system (X, f), and its relation with Ultrafilter. The concept of Profile, like Tangles and Weak tangles [49], is deeply connected with various graph parameters and has been extensively studied (e.g. [28, 50, 51]). We present the definition below.

**Definition 13 (cf. [49]):** A set $P$ of separations is consistent of which if $(C, D) \leq (A, B) \in P$ implies $(D, C) \notin P$. Note that this does not imply $(C, D) \in P$. It may also happen that $P$ contains neither $(C, D)$ nor $(D, C)$.

**Definition 14. (cf. [49]):** Let $X$ represent a finite set and $f$ denote a symmetric submodular function delineated over $X$. A Set $P$ of separations of a connectivity system $(X, f)$ is a profile of order $k + 1$ if it satisfies:

(P0) For all $(A, B) \in P$, $f(A) \leq k$,

(P1) For every separation $(A, B)$ of $X$ of order $\leq k$, one of $(A, B), (B, A)$ is an element of $P$,

(P2) P is consistent: So if $(A_2, B_2) \in P$, $(A_1, B_1) \leqq (A_2, B_2)$ , order of $(A_1, B_1) \leq k$ then $(A_1, B_1) \in P$,

(P3) For all $(A_1, B_1), (A_2, B_2) \in P : (A_1 \cap A_2, B_1 \cup B_2) \notin P$: So if $(A_1, B_1) \in P$, $(A_2, B_2) \in P$, order of $(A_1 \cup A_2, B_1 \cap B_2) \leq k$ then $(A_1 \cup A_2, B_1 \cap B_2) \in P$. We also consider about "non-principal " profile in this paper. A profile of order $k + 1$ is non-principal if it satisfies:

(P4) For each element $e \in X$ such that $f(\{e\}) \leq k$ , $(e, E \setminus \{e\}) \in P$. By modifying axiom (P3) in the above profile definition to the following (SP3), we can define a Linear-Profile of separations of order $k + 1$. (SP3) If $(A_1, B_1) \in P$ , $e \in X$, $f(\{e\}) \leq k$, then $(A_1 \cap (X \setminus \{e\}), B_1 \cup \{e\}) \notin P$::

Let us now articulate the relationship between the aforementioned Ultrafilter and Tangle. The following theorem clearly holds true.

**Theorem 15:** Let $X$ be a finite set and let $f$ be a symmetric submodular function defined on $X$. The existence of the following are equivalent conditions:

- A non-principal profile of separations in a connectivity system $(X, f)$ of order $k + 1$.

- A tangle of separations in a connectivity system *(X, f)* of order $k + 1$.
- An ultrafilter of separations in a connectivity system *(X, f)* of order $k + 1$..

**Theorem 16:** Let $X$ be a finite set and let $f$ be a symmetric submodular function defined on $X$. The existence of the following are equivalent conditions:
- A non-principal linear profile of separations in a connectivity system *(X, f)* of order $k + 1$.
- A linear tangle of separations in a connectivity system *(X, f)* of order $k + 1$.
- A single ultrafilter of separations in a connectivity system *(X, f)* of order $k + 1$.

## 6. Conclusion and Future Tasks

As for our future endeavors, we would like to outline the following tasks:
- Continuing to explore proof methods for open problems can be identified as one of the future challenges that need to be addressed.
- We plan to investigate various graph parameters and their related concepts using separations and partitions, while examining the use of the submodular partition function [29, 30] to reevaluate ultrafilters and tangles.
- We plan to reconsider loose tangle [26, 30] using Separation and Partition.
- We plan to reconsider Quasi-Ultrafilter [56-58] using Separation and Partition.
- As previously mentioned, our ultimate goal is to convincingly demonstrate that a weak ultrafilter exists if and only if a non-principal bramble exists.
- We reconsider relationship between Weak Tangle [49] and Ultrafilter.
- We consider ultraproducts in a connectivity system. To put it simply, an ultraproduct extracts properties that hold "almost everywhere." We aim to investigate the relationship between these concepts and graph width parameters.
- We will consider about filter base of separations of a connectivity system *(X, f*. We present the definition below. In simple terms, for any two elements in a filterbase, there is a "smaller" element also in the filterbase (cf. [52, 53]).

**Definition 17:** A filter base $B$ of separations of a connectivity system *(X, f) of order* $k + 1$ is satisfies:
(FB1) A filter base $B$ is not an empty set.
(FB2) If $(A_1, B_1), (A_2, B_2) \in B$, order of $(A_3, B_3) \leq k$, there exists $(A_3, B_3) \in B$ such that $(A_3, B_3) \leq (A_1, B_1)$ and $(A_3, B_3) \leq (A_2, B_2)$.

## Acknowledgments

I humbly express my sincere gratitude to all those who have extended their invaluable support, enabling me to successfully accomplish this paper.